\documentclass{article}
\usepackage{amssymb}
\usepackage{amsmath, graphicx}

\begin{document}

\begin{center}
\textbf{On moments of Pitman estimators: the case of fractional Brownian
motion}
\begin{equation*}
\end{equation*}%
Alexander Novikov \footnote{%
University of Technology, Sydney and Steklov Mathematical Institute of RAS,
Moscow. Present address: PO Box 123, Broadway, School of Mathematical
Sciences, University of Technology, Sydney, NSW 2007, Australia;
e-mail:Alex.Novikov@uts.edu.au},~Nino Kordzakhia \footnote{%
Department of Statistics, Macquarie University, North Ryde, NSW, 2109,
Australia;\; \; e-mail:Nino.Kordzakhia@mqu.edu.au} and Timothy Ling \footnote{%
University of Technology, Sydney, NSW, 2007, Australia; e-mail:Timothy.Ling@uts.edu.au}
\end{center}

\begin{equation*}
\end{equation*}%
\textbf{Abstract. }In some non-regular statistical estimation problems, the
limiting likelihood processes are functionals of fractional Brownian motion
(fBm) with Hurst's parameter $H,~0<H\leq 1$.~In this paper we present
several analytical and numerical results on the moments of Pitman
estimators represented in the form of integral functionals of fBm. %
We also provide Monte Carlo  simulation results for variances of Pitman and asymptotic maximum likelihood 
estimators.

\begin{center}
\textbf{1. Introduction}
\end{center}

Pitman estimators  (\cite{Pit}), also known as Bayesian estimators with a constant prior
on the real line (\cite{Berg}), for parameters of stochastic processes are
optimal under various continuous- and discrete-time settings, \cite{IH70}%
, \cite{IH81}. For example, one may consider the estimation problem of
parameter $\theta \ $ by observing the diffusion process $X=\{X_t,\   0\leq t \leq T\}$ that is a 
solution of  stochastic differential equation
\begin{equation*}
dX_{t}=s(X_{t}, t, \theta )dt+\sigma (X_{t})dW_{t},~0\leq t\leq T,
\end{equation*}%
where the drift $s(x, t,\theta )$ is a non-regular function, e.g. $%
s(x, t, \theta )=|x-\theta |^{p}, \; %
 p<\frac{1}{2},$ or $s(x, t, \theta )=I\{\theta>t \}$. For such non-regular
statistical estimation problems it is a typical situation when the
respective limit likelihood process $Z_{t}$ is generated by a fractional
Brownian motion (fBm)\ $W_{t}^{H}$ with Hurst's parameter $H\in (0,1],$
namely%
\begin{equation*}
Z_{t}=e^{W_{t}^{H}-\frac{1}{2}|t|^{2H}},~t\in R, \; R=(-\infty ,\infty ),
\end{equation*}%
see \cite{K}, \cite{F10}. Note that the case $H=\frac{1}{2}$ appears in a
study of a change point problem for a Brownian motion (Bm) in \cite{IH70}, 
\cite{IH81} and processes with a time delay in \cite{GuK}. The case $%
H\neq \frac{1}{2}$ appears in various continuous-time settings, see (Chapter
3, \cite{K}) and \cite{F10}, and discrete-time frameworks, \cite{IH81}, \cite%
{Husk}.

Distributional properties of Pitman estimators for large sample sizes have
not been studied in much detail. In this paper, in continuation of our
results from \cite{NK}, we study the limit distribution of Pitman
estimators, which can be defined as the distribution of a random variable%
\begin{equation}
\zeta _{H}=\int_{_{_{-\infty }}}^{^{^{\infty }}}tq_{t}dt,  \label{dzeta}
\end{equation}%
where%
\begin{equation}
q_{t}=Z_{t}(\int_{_{_{-\infty }}}^{^{^{\infty }}}Z_{u}du)^{-1},  \label{q(t)}
\end{equation}%
$\zeta _{H}$ represents a conditional expectation with respect to aposteriory density $q_{t}$.
Recall that $W^{H}=\{W_s^H, \;  s\in R\}$ is a Gaussian process with continuous trajectories%
\begin{equation*}
W_{0}^{H}=0,E(W_{s}^{H})=0,\,E|W_{s}^{H}-W_{t}^{H}|^{2}=|s-t|^{2H},\,s\in
R,~t\in R.
\end{equation*}%
This implies that the covariance function of $W_{s}^{H}$ is%
\begin{equation*}
R(t,s):=E(W_{t}^{H}W_{s}^{H})=\frac{1}{2}(|t|^{2H}+|s|^{2H}-|t-s|^{2H}).
\end{equation*}%
Note that, even \ when $H=\frac{1}{2}$, i.e. in the case of a standard Bm
 neither the distribution  nor the moments $%
\zeta _{\frac{1}{2}}\ $(greater than 2), \ are known in an explicit form.
For other cases, except $H=1$, the essential difficulty in studying the
functionals of fBm is due to the fact that $W^{H}$ is not a
semimartingale and therefore, standard tools of stochastic calculus (based
on the Ito formula) are not applicable.

The case $H=1$ corresponds to the regular statistical estimation problems
where the limit distribution is normal, $\zeta _{1}\sim N(0,1).$

In this paper we obtain several results on the variance and higher moments
of $\zeta _{H},~0<H\leq 1,$ using the measure transformation technique and
Gaussian property of fBm. In \cite{NK} we showed that, for
$H>0.309...$,  the random variable (r.v.) $|\zeta _{H}|^{2H}$
is \textit{exponentially bounded} i.e. there exists~a constant~$\alpha _{H}>0
$ such that 
\begin{equation}
Ee^{\delta |\zeta _{H}|^{2H}}<\infty ~\text{ }\text{\emph{for}}~~\delta
<\alpha _{H}.  \label{a18}
\end{equation}%
This result implies, of course, finiteness of all moments of $\zeta _{H}$.
In Section 2 we improved this result (see Theorem 1) by showing that (\ref{a18}) 
does hold for all $H\in (0,1].$  Improvement is achieved thanks to 
application of the measure transformation technique, see Lemma 1 in Section
2. Note that Lemma 1 also will be used in the proof of Theorem 2 in Section
3 which presents a general identity for expectations of functions of $\zeta
_{H}$. Then, using the obtained identity, we derive a useful representation
for the variance of $\zeta _{H}$ when $H\in (0,1]$,  (see Corollary 1). 
Corollary 2 provides a lower
bound for the moments$~E\zeta _{H}^{k},~k=2,4,..~.$

In Section 4, Theorem 3, for the case $H\in \lbrack \frac{1}{2},1)$ (see Theorem 3) we also
derive a new representation for $Var(\zeta _{H})$ in terms of the function%
\begin{equation}
g(m):=E\log (\int_{_{_{-\infty }}}^{^{^{\infty }}}Z_{u}e^{mu}du).
\label{g(m)}
\end{equation}%
This result was formulated  in \cite{NK} without proof. Earlier in \cite{NK}
we derived another expression for $Var(\zeta _{H})$ in terms of the
function 
\begin{equation*}
g(m_{1},m_{2})=E\log \int_{_{_{0}}}^{^{^{\infty
}}}(Z_{u}e^{-m_{1}u}+Z_{-u}e^{-m_{2}u})du.
\end{equation*}%
In  \cite{NK} was shown that for the case $H=\frac{1}{2}~$ the function $%
g(m_{1},m_{2})$ (as well as  $g(m)=g(-m,m)$) can be
expressed in terms of the PolyGamma function leading to much shorter
derivation of the following result from \cite{RS}%
\begin{equation}
Var(\zeta _{\frac{1}{2}})=16Zeta[3]\thickapprox 19.23,  \label{v0}
\end{equation}%
where $Zeta[k]$ is the Riemann-zeta function, see details in \cite{NK}.

In Section 5 we present Monte Carlo simulation results for $Var(\zeta _{H})$ and the variance 
of asymptotic maximum likelihood estimator which is the {\it{argmax}} of $W^{H}=\{W^H_s, s\in R\}.$

\begin{center}
\textbf{2. Exponentially boundedness of }$|\zeta _{H}|^{2H}$
\end{center}

The following lemma on a measure transformation for Gaussian processes (and
hence for fBm $W_{t}^{H},\,t\in R,~$as well) plays key role in the proof
of Theorem 1 below.

We formulate Lemma 1 in terms of a Gaussian system $(\xi
,\{X_{s}\},s\in D)$ (see \cite{Shir}) defined on probability space ($\Omega
,F,P)$. Recall that it means that $(\xi ,\{X_{t_{i}}\},t_{i}\in D,i=1,...,n)$
is a Gaussian vector for any $n$. We use the upper index, e.g., $Q,$ to indicate that
expectations are taken with respect to a measure $Q,~$so $E^{Q}(.)$ is used
for the expectation with respect to measure\emph{\ }$Q$. Lemma 1 gives a new result which belongs 
to a group of results broadly known as Cameron-Martin-Girsanov-Maruyama-... type measure transformations.

\textbf{Lemma 1.}\emph{~Let }$(\xi ,\{X_{s}\},s\in D)$\emph{\ be a Gaussian
system on a probability space (}$\Omega ,F,P)$\emph{. Set }%
\begin{equation*}
E^{P}(\xi )=0,~Var^{P}(\xi )=\sigma ^{2}
\end{equation*}%
\emph{and consider the measure transformation}%
\begin{equation*}
Q(A)=E^{P}I(A)e^{\xi -\frac{\sigma ^{2}}{2}}.
\end{equation*}%
\emph{\ Then on the probability space (}$\Omega ,F,Q):$

\emph{1)\ the system }$(\xi ,\{X_{s}\},s\in D)$\emph{\ is Gaussian;}

\emph{2)}%
\begin{equation*}
E^{Q}X_{s}=E^{P}X_{s}+Cov^{P}(\xi
,X_{s}),~Cov^{Q}(X_{t},X_{s})=Cov^{P}(X_{t},X_{s}).
\end{equation*}%
\textbf{Proof.} Property 1) is a consequence of the definition of a Gaussian
system and the fact that any linear transformation of a Gaussian vector is a
Gaussian vector.

To check the second property one should write out the joint moment
generating function of  $X_{s}~$and$~X_{t}$ with respect to measure $Q$
for $t,s\in R$%
\begin{equation*}
E^{Q}e^{z_{1}X_{t}+z_{2}X_{s}}=
\end{equation*}%
\begin{equation}
=\exp \{z_{1}E^{P}X_{s}+z_{2}E^{P}X_{t}-\frac{\sigma ^{2}}{2}+\frac{1}{2}%
Var^{P}(\xi +z_{1}X_{s}+z_{2}X_{t})\}.  \label{mgf1}
\end{equation}%
Since%
\begin{equation*}
\frac{1}{2}Var^{P}(\xi +z_{1}X_{s}+z_{2}X_{t})=\frac{1}{2}(\sigma
^{2}+z_{1}^{2}Var^{P}(X_{s})+z_{2}^{2}Var^{P}(X_{t}))+
\end{equation*}%
\begin{equation*}
+z_{1}Cov^{P}(\xi ,X_{s})+z_{2}Cov^{P}(\xi
,X_{t})+z_{1}z_{2}Cov^{P}(X_{t},X_{s}),
\end{equation*}%
differentiating in (\ref{mgf1})~with respect to $z_{1}$ and $z_{2}$ we obtain%
\begin{eqnarray*}
E^{Q}X_{s} &=&\frac{\partial }{\partial z_{1}}%
E^{Q}e^{z_{1}X_{s}+z_{2}X_{t}}| _{\substack{  \\ z_{1}=z_{2}=0}}%
=E^{P}X_{s}+Cov^{P}(\xi ,X_{s}), \\
E^{Q}X_{t}X_{s} &=&\frac{\partial ^{2}}{\partial z_{1}\partial z_{2}}%
E^{Q}e^{z_{1}X_{t}+z_{2}X_{s}}|_{\substack{  \\ z_{1}=z_{2}=0}} \\
&=&(E^{P}X_{s}+Cov^{P}(\xi ,X_{s}))(E^{P}X_{t}+Cov^{P}(\xi
,X_{t}))+Cov^{P}(X_{t},X_{s})= \\
&&E^{Q}X_{t}E^{Q}X_{s}+Cov^{P}(X_{t},X_{s}).
\end{eqnarray*}
This completes the proof.

While dealing with integrals of a Gaussian process $X_{s}$ we assumed
that there is a progressively measurable modification of $X_{s}$ such that
 integrals are well defined.

For proving Theorem 1 we also need the following lemma.

\textbf{Lemma 2. }\emph{Let }$X_{s}$\emph{\ be a Gaussian process with }$%
EX_{s}=0$\emph{. Then for any }$t>0$\emph{\ and }$r>0$%
\begin{equation}
E(\int_{0}^{t}e^{X_{s}}ds)^{-r}\leq \frac{1}{t^{r}}\exp \{\frac{r^{2}}{2t^{2}%
}Var(\int_{0}^{t}X_{s}ds)\}.  \label{Molch}
\end{equation}
\textbf{Proof.} Applying Jensen's inequality we obtain for any $t>0$%
\begin{equation*}
\int_{0}^{t}e^{X_{s}}ds\geq t\exp \{\frac{1}{t}\int_{0}^{t}X_{s}ds\},
\end{equation*}%
and hence for any $t>0$ and $r>0$ 
\begin{equation*}
E(\int_{0}^{t}e^{X_{s}}ds)^{-r}\leq \frac{1}{t^{r}}E\exp \{-\frac{r}{t}%
\int_{0}^{t}X_{s}ds\}=\frac{1}{t^{r}}\exp \{\frac{r^{2}}{2t^{2}}%
Var(\int_{0}^{t}X_{s}ds)\}.
\end{equation*}

\textbf{Remark 1.} Integrals  of  type $ E(\int_{0}^{t}e^{X_{s}}ds)^{-1}$  were considered in \cite{KawaTanaka}, \cite{Molch} and \cite{Kawa1}.

\textbf{Theorem 1. }\emph{For any}\textbf{\ }$H\in (0,1],~$\emph{there exists a
positive number }$\alpha _{H}$\emph{\ such that} 
\begin{equation}
Ee^{\delta |\zeta _{H}|^{2H}}<\infty ~~\emph{for\ }\delta <\alpha _{H}.
\label{alpha}
\end{equation}
\textbf{Proof. }If\textbf{\ }$H=1$ then it is $\alpha _{1}=\frac{1}{2}$ in (%
\ref{alpha}) because it is well known that in this case $\zeta _{1}\sim
N(0,1)$ and therefore 
\begin{equation}
Ee^{\delta \zeta _{1}^{2}}<\infty ~~\emph{for}~~\delta <\frac{1}{2};~Ee^{%
\frac{1}{2}\zeta _{1}^{2}}=\infty .  \label{H=1}
\end{equation}
In the case $H\in (0,1)$ we need to find a proper estimate for the expectation of $Eq_{t}$
which leads to (\ref{alpha}).

Note that the function $e^{\delta |x|^{2H}}$ is a convex function for  $H\geq \frac{1}{2}$, \ and 
for  $H\in (0,\frac{1}{2})$ it is dominated by convex function $\max (C _{\delta ,H},e^{\delta |x|^{2H}})$, 
where $C_{\delta ,H}$ is a sufficiently large number.
The random process $q_{t}$ represents
a density function as it is a normalised nonnegative function of $t$ 
such that $~\int_{-\infty }^{\infty
}q_{t}dt=1.$ Thus by Jensen's inequality,  for any $\delta \geq 0$,  from (\ref%
{dzeta}) we have 
\begin{eqnarray*}
e^{\delta |\zeta _{H}|^{2H}} &\leq &C_{\delta ,H}+\int_{-\infty }^{\infty
}e^{\delta |t|^{2H}}q_{t}dt \\
&\leq &C_{\delta ,H}+\int_{-1}^{1}e^{\delta |t|^{2H}}q_{t}dt+\int_{-\infty
}^{\infty }I\{|t|>1\}e^{\delta |t|^{2H}}q_{t}dt\leq
\end{eqnarray*}%
\begin{equation*}
\leq C _{\delta ,H}+e^{\delta }+\int_{-\infty }^{\infty }I\{|t|>1\}e^{\delta
|t|^{2H}}q_{t}dt,
\end{equation*}%
where the constant $ C _{\delta ,H} \geq 0,$ and $ C _{\delta ,H}=0$ in the case $H\geq 
\frac{1}{2}$; above we also used the fact that $\int_{-\infty }^{\infty
}q_{t}dt=1$. In view of the symmetry property of fBm in distributional sense
\begin{equation}
\{W_{u}^{H},u\geq 0\}\overset{d}{=}\{W_{-u}^{H},u\geq 0\},
\label{trin}
\end{equation}%
we have%
\begin{equation*}
Ee^{\delta |\zeta _{H}|^{2H}}\leq c_{\delta }+e^{\delta }+2\int_{1}^{\infty
}e^{\delta t^{2H}}Eq_{t}dt.
\end{equation*}%
For finding a proper upper bound for~$Eq_{t}$ for $t\geq 1$ we use Lemma 1
with%
\begin{equation*}
\xi =\lambda W_{t}^{H},~X_{s}=W_{s}^{H},
\end{equation*}%
where $\lambda $ is a real number. Then%
\begin{equation*}
E^{P}(X_{s})=0,\sigma ^{2}=\frac{\lambda ^{2}t^{2H}}{2};~E^{Q}(W_{t}^{H})=%
\lambda R(t,s),~R^{Q}(t,s)=R(t,s).
\end{equation*}%
This means that with respect to the measure $Q$ the process\newline
$\{W_{t}^{H}-\lambda R(t,s),t\in R\}$ is a (standard) fBm. Using this fact
we obtain%
\begin{equation*}
Eq_{t}=Ee^{\lambda W_{t}^{H}-\frac{\lambda ^{2}}{2}t^{2H}}e^{(1-\lambda
)W_{t}^{H}-\frac{(1-\lambda ^{2})}{2}t^{2H}}(\int_{-\infty }^{\infty
}e^{W_{s}^{H}-\frac{|s|^{2H}}{2}}ds)^{-1}=
\end{equation*}%
\begin{equation*}
E^{Q}e^{(1-\lambda )W_{t}^{H}-\frac{(1-\lambda ^{2})}{2}t^{2H}}(\int_{-%
\infty }^{\infty }e^{W_{s}^{H}-\frac{|s|^{2H}}{2}}ds)^{-1}.
\end{equation*}%
Applying Lemma 1 we have%
\begin{equation*}
Eq_{t}=Ee^{(1-\lambda )W_{t}^{H}+\lambda (1-\lambda )t^{2H}-\frac{(1-\lambda
^{2})}{2}t^{2H}}(\int_{-\infty }^{\infty }e^{W_{s}^{H}+\lambda R(t,s)-\frac{%
|s|^{2H}}{2}}ds)^{-1}=
\end{equation*}%
\begin{equation}
e^{-\frac{(1-\lambda )^{2}}{2}t^{2H}}Ee^{(1-\lambda
)W_{t}^{H}}(\int_{-\infty }^{\infty }e^{W_{s}^{H}+t^{2H}f(s/t)}ds)^{-1},
\label{aaa}
\end{equation}%
where 
\begin{equation*}
t^{2H}f(s/t):=\lambda R(t,s)-\frac{|s|^{2H}}{2}.
\end{equation*}%
Note that 
\begin{equation}
f(u)=\frac{\lambda +(\lambda -1)u^{2H}-\lambda |1-u|^{2H}}{2}%
,~f(0)=0,~f(1)=\lambda -1/2.
 \label{qqqq}
\end{equation}%
To simplify the exposition of the proof we choose $\lambda =\frac{1}{2}~$%
(although it seems that somewhat  better estimator for $\alpha _{H}~$%
can be obtained with a proper choice of $\lambda ~$depending on $H$). Then assuming
 $\lambda =\frac{1}{2}~$\ we can rewrite (\ref{aaa}) as follows%
\begin{equation}
Eq_{t}=e^{-\frac{1}{8}t^{2H}}Ee^{\frac{1}{2}W_{t}^{H}}(\int_{-\infty
}^{\infty }e^{W_{s}^{H}+t^{2H}f(s/t)}ds)^{-1},  \label{aaa1}
\end{equation}%
where 
\begin{equation}
f(u)=f(1-u)=\frac{1-|u|^{2H}-|1-u|^{2H}}{4},~f(0)=0,~f(1)=0.  \label{fuf1}
\end{equation}%
It is easy to see that for $H\in (0,1/2)$ the function~$f(u),0<u<1$ is
negative; for $H\in (1/2,1)$ the function~$f(u),0<u<1$ is positive.

Applying the inequality$\ e^{x/2}\leq (1+e^{x})/2$ for the term $%
e^{W_{t}^{H}/2}$ in (\ref{aaa1}) and reducing the range of integration to $%
s\in \lbrack 0,t]$ instead of $s\in R$ for the integral we obtain 
\begin{equation*}
Eq_{t}\leq \frac{1}{2}e^{-\frac{1}{8}t^{2H}}[E(%
\int_{0}^{t}e^{W_{s}^{H}+t^{2H}f(s/t)}ds)^{-1}+E(%
\int_{0}^{t}e^{W_{s}^{H}-W_{t}^{H}+t^{2H}f(s/t)}ds)^{-1}].
\end{equation*}%
In view of following translation-invariance property of fBm%
\begin{equation}
\{W_{s}^{H}-W_{t}^{H},s\in R\}\overset{d}{=}\{W_{t-s}^{H},~s\in R\},
\label{trin}
\end{equation}%
the expectation $E(\int_{0}^{t}e^{W_{s}^{H}-W_{t}^{H}+t^{2H}f(s/t)}ds)^{-1}$
 is equal to $E(%
\int_{0}^{t}e^{W_{s}^{H}+t^{2H}f(1-s/t)}ds)^{-1}$  after change of variable $t-s$  to $s$. Thus we get the estimate%
\begin{equation*}
Eq_{t}\leq \frac{1}{2}e^{-\frac{1}{8}t^{2H}}[E(%
\int_{0}^{t}e^{W_{s}^{H}+t^{2H}f(s/t)}ds)^{-1}+E(%
\int_{0}^{t}e^{W_{s}^{H}+t^{2H}f(1-s/t)}ds)^{-1}].
\end{equation*}%
Since $f(u)=f(1-u)$ the above inequality can be now rewritten as follows 
\begin{equation}
Eq_{t}\leq e^{-\frac{1}{8}t^{2H}}E(%
\int_{0}^{t}e^{W_{s}^{H}+t^{2H}f(s/t)}ds)^{-1}.  \label{ll}
\end{equation}%
\textbf{Consider now the case} $H\in \lbrack 1/2,1).~$Then it is easy to see 
$f(s/t)\geq 0$ and therefore by Lemma 2 for $t\geq 1$%
\begin{equation*}
Eq_{t}\leq e^{-\frac{1}{8}t^{2H}}E(\int_{0}^{1}e^{W_{s}^{H}}ds)^{-1}<\infty .
\end{equation*}%
Combining all estimates obtained above for
 $H\in \lbrack \frac{1}{2},1)$ and $\delta <\frac{1}{8}$ we have 
\begin{equation*}
Ee^{\delta |\zeta _{H}|^{2H}}\leq c_{\delta }+e^{\delta
}+2E(\int_{0}^{1}e^{W_{s}^{H}}ds)^{-1}\int_{1}^{\infty }e^{-(\frac{1}{8}%
-\delta )t^{2H}}dt<\infty .
\end{equation*}%
Since $E(\int_{0}^{1}e^{W_{s}^{H}}ds)^{-1}$ is finite due to the result of
Lemma 2 we have proved Theorem 1 with $\alpha _{H}\geq \frac{1}{8}.$

\textbf{The case }$H\in (0,1/2)$.

Obviously, the function $f(u)$ (defined in (\ref{fuf1})) is decreasing on
the interval $u\in (0,\frac{1}{2})$. This fact implies that $f(s/t)\geq
f(\varepsilon )$ for all $s\in (0,\varepsilon t)~$and any $\varepsilon \in
(0,\frac{1}{2})$, where $f(\varepsilon )=-(\varepsilon ^{2H}/4)(1+o(1))~$\ as $%
\varepsilon \rightarrow 0$.~Thus, from (\ref{ll}) we have 
\begin{equation*}
Eq_{t}\leq e^{-(\frac{1}{8}+f(\varepsilon ))t^{2H}}E(\int_{0}^{t\varepsilon
}e^{W_{s}^{H}}ds)^{-1}.
\end{equation*}%
Again,  the fact that $E(\int_{0}^{t\varepsilon }e^{W_{s}^{H}}ds)^{-1}$ is a bounded
and decreasing function of $t\in \lbrack 1,\infty )$ easily follows from Lemma 2.

Combining all estimates obtained above for
$H\in (0,\frac{1}{2}]$ we have the following estimate 
\begin{equation*}
Ee^{\delta |\zeta _{H}|^{2H}}\leq c_{\delta }+e^{\delta
}+2E(\int_{0}^{\varepsilon }e^{W_{s}^{H}}ds)^{-1}\int_{1}^{\infty }e^{-(%
\frac{1}{8}+f(\varepsilon )-\delta )t^{2H}}dt
\end{equation*}%
where the right-hand side is obviously finite when $\varepsilon ~$\ is small
enough and $\delta <\frac{1}{8}.$

This completes the proof of Theorem 1 with $\alpha _{H}\geq \frac{1}{8}$ for
all $H\in (0,1].$

\textbf{Remark 2.}\ Using a different approach in \cite{NK} we
found that
\begin{equation*}
\alpha _{H}\geq \frac{4H^{2}+2H-1}{2(2H+2)(2H+1)}~~\text{\emph{for}}~H>(%
\sqrt{5}-1)/4=0.3090...~~.
\end{equation*}%
Comparing this estimate with that one obtained in the proof of Theorem 1 we
get the following  lower bounds:%
\begin{equation*}
\alpha _{H}\geq \frac{1}{8}~~\text{\emph{for}~}H\in (0,H_{0}],~\alpha
_{H}\geq \frac{4H^{2}+2H-1}{2(2H+2)(2H+1)}~\text{\emph{for}~}H\in \lbrack
H_{0},1),
\end{equation*}%
where $H_{0}=(\sqrt{73}-1)/12=0.6287...,~H_{0}$ is the largest root of the equation 
\begin{equation*}
\frac{4H^{2}+2H-1}{2(2H+2)(2H+1)}=\frac{1}{8}.
\end{equation*}
Recall that for $H=1$ the index $\alpha _{1}=\frac{1}{2}$,  see (\ref{H=1}).

\textbf{Conjecture.} \emph{There exists an index }$\alpha _{H}$\emph{\ such
that}%
\begin{equation*}
Ee^{\delta |\zeta _{H}|^{2H}}<\infty ~~\text{for }\delta <\alpha
_{H},~Ee^{\delta |\zeta _{H}|^{2H}}=\infty \text{ for }\delta >\alpha _{H}.
\end{equation*}

\textbf{Remark 3. }This conjecture is motivated by the result of Theorem 1
and the following result on the limit distribution of Maximum Likelihood
Estimator (MLE) $\xi_H$ for the case $H=1/2$.

It is well
known that the distribution of   $\xi_{\frac{1 }{2}}$  coincides with the distribution of a location of maximum of two-sided
Brownian motion and is
\begin{equation}
P(|\xi_{\frac{1}{2}}|>t)=(t+5)\Phi (-\frac{\sqrt{t}}{2})-\sqrt{\frac{2t}{\pi }}e^{-%
\frac{t}{8}}-3e^{t}\Phi (-\frac{3\sqrt{t}}{2}),  \label{yao}
\end{equation}%
where $\Phi (t)$ is a standard normal distribution. This result
can be easily derived from the papers \cite{Shepp} and \cite{Yao}. Using the
well known formula $\Phi (-x)=\sqrt{\frac{1}{2\pi x^{2}}}e^{-\frac{x^{2}}{2}%
}(1+o(1)),\  x\rightarrow \infty ,$  we have
\begin{equation*}
P(|\xi_{\frac{1}{2}}|>t)=\sqrt{\frac{32}{\pi t}}e^{-\frac{t}{8}}(1+o(1)),\,\,t%
\rightarrow \infty 
\end{equation*}%
and, hence, 
\begin{equation*}
Ee^{\delta |\xi_{\frac{1}{2}}|}<\infty ~~\text{for }\delta <\frac{1}{8},~Ee^{\frac{1%
}{8}|\xi_{\frac{1}{2}}|}=\infty .
\end{equation*}
Note that,  from  (\ref{yao}) one can directly obtain  
\begin{equation}
Var(\xi _{\frac{1}{2}})=26.  \label{Tere26}
\end{equation}%
This result appeared in \cite{Tere} for the first time.

\textbf{Remark 4.} The reviewer of this paper indicated that the existence
of exponential moments for ${|\zeta _{H}|}^{2H}$ can be retraced from the general
results of  Ibragimov-Hasminski theory (see \cite{IH81}) in combination with some results
from \cite{Dach}.

\begin{center}
\textbf{3. Identities for expectations of functions of  $\zeta _{H}$.}
\end{center}

\textbf{Theorem 2. }\emph{Let }$G(\zeta _{H})$\emph{\ be a mesurable bounded
function of }$\zeta _{H},~H\in (0,1].~$\emph{Then}%
\begin{equation}
EG(\zeta _{H})=\int_{_{_{-\infty }}}^{^{^{\infty }}}EG(\zeta _{H}-t)q_{t}dt.
\label{Ruhin}
\end{equation}
\textbf{Proof. }Using (\ref{dzeta}) and Lemma 1 with%
\begin{equation*}
\xi =W_{t}^{H},~X_{s}=W_{s}^{H},
\end{equation*}%
we obtain%

\begin{equation*}
\int_{_{_{-\infty }}}^{^{^{\infty }}}EG(\zeta
_{H}-t)q_{t}dt=\int_{_{_{-\infty }}}^{^{^{\infty }}}E^{Q}G(\zeta
_{H}-t)(\int_{_{_{-\infty }}}^{^{^{\infty }}}Z_{s}ds)^{-1}dt
\end{equation*}%

(using Lemma 1 after symplifying we get)

\begin{equation*}
\begin{split}
=\int_{_{_{-\infty }}}^{^{^{\infty }}}EG\left(\int_{_{_{-\infty }}}^{^{^{\infty
}}}se^{W_{s}^{H}-\frac{1}{2}|t-s|^{2H}}ds\left(\int_{_{_{-\infty }}}^{^{^{\infty
}}}e^{W_{s}^{H}+\frac{1}{2}|t|^{2H}-\frac{1}{2}|t-s|^{2H}}ds\right)^{-1}-t\right)  \\
\times \left(%
\int_{_{_{-\infty }}}^{^{^{\infty }}}e^{W_{s}^{H}+\frac{1}{2}|t|^{2H}-\frac{1%
}{2}|t-s|^{2H}}ds\right)^{-1}dt.
\end{split}
\end{equation*}%
Next, using (\ref{trin}) we have%

\begin{equation*}
\int_{_{_{-\infty }}}^{^{^{\infty }}}EG(\zeta _{H}-t)q_{t}dt=
\end{equation*}%
\begin{equation*}
\begin{split}
\int_{_{_{-\infty }}}^{^{^{\infty }}}EG\left(\int_{_{_{-\infty }}}^{^{^{\infty
}}}se^{W_{s-t}^{H}-\frac{1}{2}|t-s|^{2H}}ds\left(\int_{_{_{-\infty
}}}^{^{^{\infty }}}e^{W_{s-t}^{H}-\frac{1}{2}%
|t-s|^{2H}}ds\right)^{-1}-t\right)\\
\times  e^{W_{t}^{H}-\frac{1}{2}|t|^{2H}}\left(\int_{_{_{-\infty
}}}^{^{^{\infty }}}e^{W_{s-t}^{H}-\frac{1}{2}|t-s|^{2H}}ds\right)^{-1}dt
\end{split}
\end{equation*}%
(setting $s-t=u)$%
\begin{equation*}
\begin{split}
=\int_{_{_{-\infty }}}^{^{^{\infty }}}EG\left(\int_{_{_{-\infty }}}^{^{^{\infty
}}}ue^{W_{u}^{H}-\frac{1}{2}|u|^{2H}}du\left (\int_{_{_{-\infty }}}^{^{^{\infty
}}}e^{W_{u}^{H}-\frac{1}{2}|u|^{2H}}du\right)^{-1}\right )\\
\times e^{W_{t}^{H}-\frac{1}{2}%
|t|^{2H}}\left(\int_{_{_{-\infty }}}^{^{^{\infty }}}e^{W_{u}^{H}-\frac{1}{2}%
|u|^{2H}}du\right)^{-1}dt\\
=\int_{_{_{-\infty }}}^{^{^{\infty }}}EG(\zeta _{H})q_{t}dt=EG(\zeta
_{H})\int_{_{_{-\infty }}}^{^{^{\infty }}}q_{t}dt=EG(\zeta _{H}).
\end{split}
\end{equation*}%
This completes the proof.

\textbf{Remark 5. }\emph{A discrete-time analog of }(\ref{Ruhin})\emph{\ for
independent identically distributed (iid) observations can be found in} (%
\cite{Bor}, Lemma 2.18.1), \cite{G}.

Further we use notations: 
\begin{equation*}
B_{i}=\int_{_{_{-\infty }}}^{^{^{\infty }}}t^{i}Z_{t}dt,~i=0,1,2...
\end{equation*}%
\begin{equation*}
A_{p}=\int_{_{_{-\infty }}}^{^{^{\infty }}}|t|^{p}q_{t}dt,~p>0.
\end{equation*}%
Due to the fact that $q_{t}=Z_{t}(\int_{_{_{-\infty }}}^{^{^{\infty
}}}Z_{u}du)^{-1}=\frac{Z_{t}}{B_{0}}$ is a density function and due to the
Holder inequality we obtain from (\ref{dzeta}) that for any $p\geq 1$ 
\begin{equation}
|\zeta _{H}|^{p}\leq A_{p}.  \label{upbound}
\end{equation}
\textbf{Corollary 1. }\emph{For any }$H\in (0,1]$ 
\begin{equation}
Var(\zeta _{H})=E\zeta _{H}^{2}=\frac{1}{2}EA_{2}.  \label{Gol}
\end{equation}
\textbf{Proof. } Let $G(\zeta _{H})=\min (|\zeta
_{H}|^{2},K)$ with a finite parameter $K>0.$ Then in view of Theorem 1 and
passing to the limit as $K\rightarrow \infty $ (using the Lebesgue theorem
and Fatou's lemma) we obtain%
\begin{equation*}
E\zeta _{H}^{2}=\int_{_{_{-\infty }}}^{^{^{\infty }}}E(\zeta
_{H}-t)^{2}q_{t}dt.
\end{equation*}%
Note that by (\ref{upbound}) for any $H\in (0,1]$%
\begin{equation*}
E\zeta _{H}^{2}=EA_{1}^{2}<\infty ,~\int_{_{_{-\infty }}}^{^{^{\infty
}}}t^{2}Eq_{t}du<\infty .
\end{equation*}%
Expanding $(\zeta _{H}-t)^{2}=\zeta _{H}^{2}-2\zeta _{H}t+t^{2}$~this implies%
\begin{eqnarray*}
E\zeta _{H}^{2} &=&\int_{_{_{-\infty }}}^{^{^{\infty }}}E(\zeta
_{H}^{2}-2\zeta _{H}t+t^{2})\frac{Z_{t}}{B_{0}}dt=E\zeta _{H}^{2}-2E\zeta
_{H}A_{1}+EA_{2} \\
&=&E\zeta _{H}^{2}-2E\zeta _{H}^{2}+EA_{2}.
\end{eqnarray*}%
After simplifying we get (\ref{Gol}).

\textbf{Remark 6. } Originally the identity (\ref{Gol}) was proved in \cite{G}
for $H > \frac{1}{2}.~$The method used in \cite{G} was based on the fact
that$~$a similar identity is valid for Pitman estimators of a location
parameter for independent identically distributed observations.

Theorem 1 can be used for derivation of various useful properties of the
distribution of $\zeta _{H}.~$As another example we present the following
result.

\textbf{Corollary 2. }\emph{For any }$H\in (0,1]~$\emph{\ and }$k=2,4,6,...~ 
$\emph{there exist constants }$c_{k}>0$\emph{\ such that}%
\begin{equation}
E\zeta _{H}^{k}\geq c_{k}^{k}EA_{k},  \label{ckey}
\end{equation}%
\emph{where }$c_{k}$\emph{\ is the unique positive root of the equation}%
\begin{equation}
(x+1)^{k}-(x-1)^{k}+2k(x^{k}-x^{k-1})=2.  \label{ckey2}
\end{equation}
\textbf{Proof. }The validity of the result for $k=2~$\ can be seen from
Corollary 1.

For the case $k\geq 4$ we apply Theorem 1 with the polynomial $G(x)=x^{k}~.$
Then we obtain 
\begin{equation*}
E\zeta _{H}^{k}=E\zeta _{H}^{k}-kE\zeta
_{H}^{k-1}A_{1}+\sum_{i=2}^{k-1}(-1)^{i}C_{k}^{i}E\zeta _{H}^{k-i}A_{i}+EA_{k},
\end{equation*}%
where $C_{k}^{i}$ are binomial coefficients. This implies 
\begin{equation*}
kE\zeta _{H}^{k}=\sum_{i=2}^{k-1}(-1)^{i}C_{k}^{i}E\zeta _{H}^{k-i}A_{i}+EA_{k}
\end{equation*}%
and hence 
\begin{equation*}
kE(\zeta _{H}^{k})\geq -\sum_{i=odd\geq 3}^{k-1}C_{k}^{i}E|\zeta
_{H}^{k-i}|A_{i}+EA_{k}.
\end{equation*}%
By the Holder inequality%
\begin{equation*}
E|\zeta _{H}^{k-i}|~A_{i}\leq (E\zeta
_{H}^{k})^{1-i/k}(EA_{i}^{k/i})^{i/k}\leq (E\zeta
_{H}^{k})^{1-i/k}(EA_{k})^{i/k},
\end{equation*}%
we have%
\begin{equation*}
kE\zeta _{H}^{k}\geq EA_{k}-\sum_{i=odd\geq 3}^{k-1}C_{k}^{i}(E\zeta
_{H}^{k})^{1-i/k}(EA_{k})^{i/k}.
\end{equation*}%
Set $x^{k}=(E\zeta _{H}^{k})/EA_{k}$. Then 
\begin{equation*}
(E\zeta _{H}^{k})^{1-i/k}/(EA_{k})^{1-i/k}=(E\zeta
_{H}^{k}/EA_{k})^{1-i/k}=x^{k-i}
\end{equation*}%
and therefore the last inequality is equivalent to%
\begin{equation*}
kx^{k}+\sum_{i=odd\geq 3}^{k-1}C_{k}^{i}x^{k-i}\geq 1.
\end{equation*}%
We can find a short expression for $\sum_{i=odd\geq 3}^{k-1}C_{k}^{i}x^{k-i}$
using the following elementary identity 
\begin{equation*}
2\sum_{i=odd\geq 3}^{k-1}C_{k}^{i}x^{k-i}=(x+1)^{k}-(x-1)^{k}-2kx^{k-1}.
\end{equation*}%
Hence we obtain%
\begin{equation*}
(x+1)^{k}-(x-1)^{k}+2k(x^{k}-x^{k-1})\geq 2.
\end{equation*}%
This implies the result.

\textbf{Remark 7. }\emph{One can easily verify that}%
\begin{equation}
c_{k}=\frac{D}{k}(1+o(1)),~k\rightarrow \infty ,  \label{cc}
\end{equation}%
\emph{where }$D=\ln(1+\sqrt{2})$\emph{\ is the unique positive root of the
equation} 
\begin{equation*}
\sinh (D)=1.
\end{equation*}%
The derivation of (\ref{cc}) is elementary and is omitted. We restrict ourselves to
illustration of accuracy of the approximation (\ref{cc}) for $k=100$,  $%
D/100\thickapprox 8.813\times 10^{-3}$,  and in this case the exact solution of (\ref{ckey2}%
)~is $c_{100}=8.841\,2\times 10^{-3}.$
\newpage
\begin{center}
\textbf{4. Representation for }$Var(\zeta _{H})$.
\end{center}

In this section for the case $H\in \lbrack \frac{1}{2},1]$ we derive another
representation for $Var(\zeta _{H})$ in terms of the function $g(m)$ defined
above in (\ref{g(m)}).

Furthermore we use the following parametrised random functions:%
\begin{equation*}
\alpha (m)=\int_{-\infty }^{\infty }ue^{mu}Z_{u}du,~\beta (m)=\int_{-\infty
}^{\infty }e^{mu}Z_{u}du~
\end{equation*}%
where $m$ is an auxiliary parameter. In these notations we have%
\begin{equation*}
\int_{-\infty }^{\infty }Z_{u}du=\beta (0),~\int_{-\infty }^{\infty
}uZ_{u}du=\alpha (0),~\alpha (m)=\frac{\partial }{\partial m}\beta (m)
\end{equation*}%
and%
\begin{equation*}
\zeta _{H}=\frac{\alpha (0)}{\beta (0)}.
\end{equation*}
Note that due to the symmetry property of fBm 
we have $g(m)=g(-m)~$ and from inequality $\log (a+b)\leq \log
(a+1)+\log (b+1),~(a>0,~b>0),~$\ we have
\begin{equation}
g(m)\leq E\log (\int_{-\infty }^{0}e^{mu}Z_{u}du+1)+E\log (\int_{0}^{\infty
}e^{mu}Z_{u}du+1).  \label{rhs}
\end{equation}
Let $m>0.$ The finiteness of the first expectation in the RHS
of (\ref{rhs}) is obvious due to the inequality $\log (x+1)\leq x$ and the
equality $EZ_{u}=1.$

The finiteness of the second expectation in the right-hand side (RHS) of (%
\ref{rhs}) for $m>0$ can be shown as follows.

Note%
\begin{eqnarray*}
\log (\int_{0}^{\infty }e^{mu}Z_{u}du+1) &=&\int_{0}^{\infty
}(\int_{0}^{s}e^{mu}Z_{u}du+1)^{-1}d\int_{0}^{s}e^{mu}Z_{u}du \\
&\leq &\int_{1}^{\infty
}e^{ms}Z_{s}(\int_{0}^{s}Z_{u}du)^{-1}ds+e^{m}\int_{0}^{1}Z_{s}(%
\int_{0}^{s}Z_{u}du)^{-1}ds \\
&\leq &\int_{1}^{\infty }e^{ms}Z_{s}(\int_{0}^{s}Z_{u}du)^{-1}ds+e^{m}.
\end{eqnarray*}%
Since $EZ_{s}(\int_{0}^{s}Z_{u}du)^{-1}=Eq_{s}\leq Ce^{-\delta s^{2H}}$ for $%
s\geq 1$ and $\delta <\frac{1}{8}$ (see the proof of Theorem 1 and (\ref{a18}%
), from now on $C$ is a generic constant) we can claim that $g(m),0<m<\frac{1%
}{8}$ is finite~(recall that we assumed $H\in \lbrack \frac{1}{2},1])$.
Obviously, $g(m)$ is a continuous function.

\textbf{Theorem 3. }\emph{Let}$~H\in \lbrack \frac{1}{2},1].$\textbf{%
\thinspace }\emph{Then the function }$g(m)$ \emph{is twice continuously
differentiable on the interval }$m\in (-1/8,1/8)$ \emph{and}%
\begin{equation*}
Var(\zeta _{H})=\frac{\partial ^{2}g(m)}{\partial m^{2}}|_{m=0}.
\end{equation*}
\textbf{Proof.} Using the notation introduced above we have 
\begin{equation*}
Var(\zeta _{H})=E\frac{\alpha ^{2}(0)}{\beta ^{2}(0)}=\lim_{m\rightarrow 0}E%
\frac{\alpha ^{2}(m)}{\beta ^{2}(m)}.
\end{equation*}%
The last equality can be justified by (\ref{a18}) and the estimate 
\begin{equation*}
E\frac{\alpha ^{2}(m)}{\beta ^{2}(m)}=E(\int_{-\infty }^{\infty }u\frac{%
e^{mu}Z_{u}}{\beta (m)}du)^{2}\leq E\int_{-\infty }^{\infty }u^{2}\frac{%
e^{mu}Z_{u}}{\beta (m)}du<\infty .
\end{equation*}%
By direct calculations we obtain for $m>0$ that%
\begin{equation*}
\frac{\partial ^{2}\log \beta (m)}{\partial m^{2}}=-\frac{\alpha ^{2}(m)}{%
\beta ^{2}(m)}+\frac{\int_{-\infty }^{\infty }u^{2}e^{mu}Z_{u}du}{\beta (m)}.
\end{equation*}%
Applying the expectation to both sides of the last equality and using
well-known theorems about differentiability of expectations with respect to
a parameter we obtain%
\begin{equation}
\frac{\partial ^{2}g(m)}{\partial m^{2}}=E\frac{\partial ^{2}\log \beta (m)}{%
\partial m^{2}}=-E\frac{\alpha ^{2}(m)}{\beta ^{2}(m)}+E\frac{\int_{-\infty
}^{\infty }u^{2}e^{mu}Z_{u}du}{\beta (m)},  \label{qq}
\end{equation}%
where the RHS is a continuous function of $m$. This implies $\frac{\partial
^{2}g(m)}{\partial m^{2}}$ is a continuous function for $m\in (0,1/8) $ and
(due to symmetry) also for $m\in (-1/8,0)$. Passing to the limit in (\ref{qq}%
) as $m\rightarrow 0$ we obtain%
\begin{eqnarray*}
\frac{\partial ^{2}g(m)}{\partial m^{2}}|_{m=0} &=&-E\frac{\alpha ^{2}(0)}{%
\beta ^{2}(0)}+E\frac{\int_{-\infty }^{\infty }u^{2}Z_{u}du}{\beta (0)} \\
&=&-Var(\zeta _{H})+2Var(\zeta _{H})=Var(\zeta _{H}).
\end{eqnarray*}
This completes the proof.

\begin{center}
\textbf{5. Modelling results. }
\end{center}
To the best of our knowledge the problem of evaluation of integral functionals numerically remains  unresolved. The only known explicit result is given by formula
(\ref{v0}). These integral functionals can be modelled using Monte-Carlo simulation method. The results of Monte-Carlo modelling for variances 
of Pitman  estimator  $\zeta _{H}$ and asymptotic MLE $\xi _{H }$ for $H\in \lbrack 0.4,1)$ are given in the Table 1. 
\newline
\begin{table}
\begin{tabular}{p{1cm}p{1.5cm}p{1.5cm}p{1.8cm}p{1.5cm}p{1.5cm}p{1.5cm}}
    \hline\noalign{\smallskip}
    $H$ & $\hat{Var}(\zeta_H)$ &  \( SE_{Var(\zeta_H)} \) 
     & $\hat{Var}(\xi_H)$ & \( SE_{Var(\xi_H)} \)\\
     \noalign{\smallskip}\noalign{\smallskip}
    0.4   & 109.682 & 0.4698 & 151.707 & 0.2145\\
    0.5   & 19.2544 & 0.0350 & 25.964   &0.0367\\
    0.6   & 6.52596 & 0.0163 & 8.63501 & 0.0386\\
    0.7   & 3.16871 & 0.0066 & 4.08858 & 0.0182\\
    0.81 & 1.82699 & 0.0032 &  2.24197 & 0.0101 \\
    0.91 & 1.28289 & 0.002   & 1.47782 & 0.0066\\
    \noalign{\smallskip}\hline\noalign{\smallskip}
    \end{tabular}
    \caption{Monte Carlo estimates for \( Var(\zeta_H) \)  and \( Var(\xi_H) \) using  \( 10^{6} \)
trajectories. Each trajectory is generated with \( 2^{18} \)
      equally-spaced discretisation points on the interval (\(- 10^{5} \), \( 10^{5}\)).}
\label{tab:zetaResults}
\end{table}
For simulation of  increments of fBm we implemented the "Circulant embedding
method" (see \cite{AnWood}) which is recognised as one of the fastest
methods for simulation of stationary Gaussian processes.

\begin{figure}
\includegraphics[height=7cm]{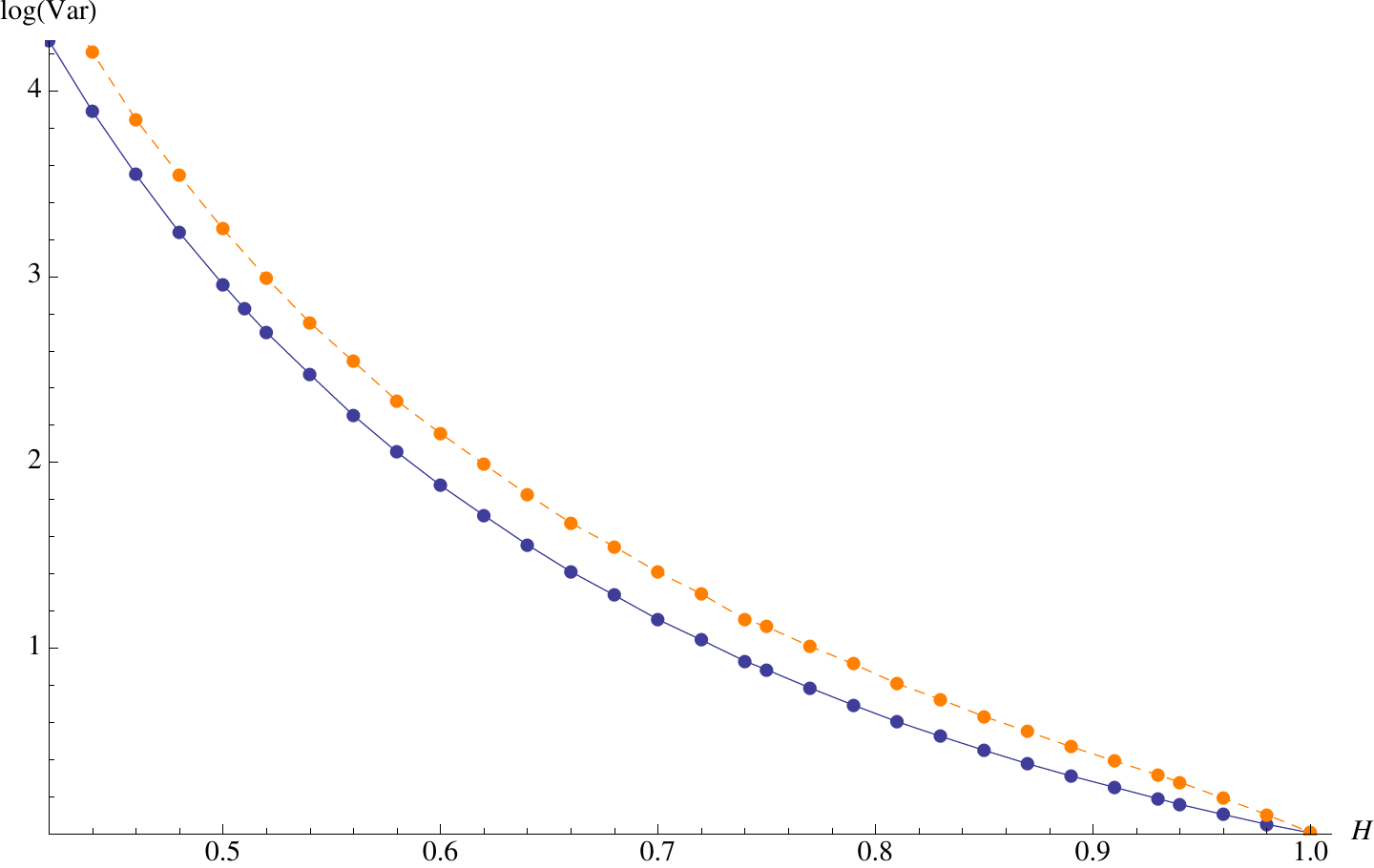}%
\caption{Values of {$Var(\zeta_H)$  (solid line) and $Var(\xi_H)$  (dashed line) for
\( H \in  [0.4,1]\)} are given on a logarithmic axis.}
\end{figure}

The graphs of $Var(\zeta _{H})$ and  $Var(\xi _{H })$ versus 
$H\in \lbrack 0.4,1)$ are plotted in Figure 1. Both $Var(\zeta _{H})$ and  $Var(\xi _{H })$ are 
monotone functions taking larger values for small values of $H$,  $Var(\zeta _{H})<Var(\xi _{H })$.
The results of  calculations agree well with formulae (\ref{v0}) and (\ref{Tere26}).
Detailed discussion of accuracy of   $Var(\zeta _{H})$  and  $Var(\xi _{H })$  are provided  in \cite{LingNov}.

\textbf{Acknowledgement.} The authors are thankful to Yuri Kutoyants for the
suggestion to study the properties of asymptotic distributions of Pitman
estimators.The  authors  also thank Julia Mishura for useful comments.

\end{document}